\documentclass[12pt]{amsart}
\usepackage{amssymb}
\usepackage{amsfonts}

\newcommand{\nc}{\newcommand}
\nc{\thref}[1]{Theorem~\ref{theo:#1}}
\nc{\selabel}[1]{\label{sect:#1}}
\nc{\seref}[1]{Section~\ref{sect:#1}}
\nc{\lelabel}[1]{\label{lemm:#1}}
\nc{\leref}[1]{Lemma~\ref{lemm:#1}}
\nc{\prlabel}[1]{\label{prop:#1}}
\nc{\prref}[1]{Proposition~\ref{prop:#1}}
\nc{\colabel}[1]{\label{coro:#1}}
\nc{\coref}[1]{Corollary~\ref{coro:#1}}
\nc{\exlabel}[1]{\label{exam:#1}}
\nc{\exref}[1]{Example~\ref{exam:#1}}
\nc{\delabel}[1]{\label{defi:#1}}
\nc{\deref}[1]{Definition~\ref{defi:#1}}
\nc{\eqlabel}[1]{\label{equa:#1}}
\nc{\relabel}[1]{\label{rema:#1}}
\nc{\reref}[1]{Lemma~\ref{rema:#1}}
\providecommand{\operatorname}[1]{\mathrm{#1}\,}
\nc{\Hom}{\operatorname{Hom}} \nc{\Mor}{\operatorname{Mor}}
\nc{\Aut}{\operatorname{Aut}} \nc{\Ann}{\operatorname{Ann}}
\nc{\Ker}{\operatorname{Ker}} \nc{\Trace}{\operatorname{Trace}}
\nc{\Char}{\operatorname{Char}} \nc{\Mod}{\operatorname{Mod}}
\nc{\End}{\operatorname{End}} \nc{\Spec}{\operatorname{Spec}}
\nc{\Span}{\operatorname{Span}} \nc{\sgn}{\operatorname{sgn}}
\nc{\Id}{\operatorname{Id}} \nc{\Com}{\operatorname{Com}}
\nc{\rank}{\operatorname{rank}}

\let\:=\colon

\newtheorem{de}{Definition}[section]
\newtheorem{lm}[de]{Lemma}
\newtheorem{pr}[de]{Proposition}
\newtheorem{co}[de]{Corollary}
\newtheorem{re}[de]{Remark}
\newtheorem{res}[de]{Remarks}
\newtheorem{te}[de]{Theorem}
\newtheorem{ex}[de]{Example}
\newtheorem{exs}[de]{Examples}

\newtheorem{prb}[de]{Problem}

\def\bex{\begin{ex}}
\def\eex{\end{ex}}
\def\bexs{\begin{exs}}
\def\eexs{\end{exs}}
\def\bl{\begin{lm}}
\def\el{\end{lm}}
\def\bc{\begin{co}}
\def\ec{\end{co}}
\def\bt{\begin{te}}
\def\et{\end{te}}
\def\bpr{\begin{pr}}
\def\epr{\end{pr}}
\def\br{\begin{re}}
\def\er{\end{re}}
\def\brs{\begin{res}}
\def\ers{\end{res}}
\def\bd{\begin{de}}
\def\ed{\end{de}}
\def\bp{\begin{prb}}
\def\ep{\end{prb}}
\def\be{\begin{equation}}
\def\ee{\end{equation}}
\def\bea{\begin{eqnarray*}}
\def\eea{\end{eqnarray*}}
\def\bp{\begin{proof}}
\def\ep{\end{proof}}

\def\qed{\hfill\Box}

\let\:=\colon
\begin{document}
\title[\hfilneg \hfil rational zeta series involving $\zeta(2n)$]
{On some rational zeta series involving $\zeta(2n)$ and binomial coefficients}

\dedicatory{}

\begin{abstract}
In this note, we give an exact formula for a general family of rational zeta series involving the coefficient $\zeta(2n)$ in terms of Hurwitz zeta values. This formula generalizes two previous formulas from a paper in \cite{Lupu-Orr}. Our method will involve derivatives polynomials for the cotangent function.  

\end{abstract}

\author{Cezar Lupu, Vlad Matei}



\maketitle

\section{Introduction}

The Riemann zeta function is defined by the absolutely convergent series

$$\displaystyle\zeta(s)=\sum_{n=1}^{\infty}\frac{1}{n^s}, \operatorname{Re}s>1.$$
\vspace{0.1cm}

In 1882, Hurwitz defined the following "shifted" zeta function,

$$\displaystyle\zeta(s; a)=\sum_{n=0}^{\infty}\frac{1}{(n+a)^s}, \operatorname{Re}s>1, 0<a\leq 1.$$

Both of them have similar properties in many aspects. For example, both of them are analytic and they have analytic continuation to the whole complex plane except for the pole $s=1$. Some particular values include $\displaystyle\zeta(-n; a)=-\frac{B_{n+1}(a)}{n+1}$, where $B_{k}(a)$ is the Bernoulli polynomial which is defined by the power series $$\displaystyle\frac{te^{xt}}{e^t-1}=\sum_{n=0}^{\infty}B_{n}(x)\frac{t^n}{n!}.$$

Also, as a special case, we have $\zeta(0; a)=\frac{1}{2}-a$. Other obvious values include $\zeta\left(s; \frac{1}{2}\right)=(2^s-1)\zeta(s)$ and $\zeta(s; a+1)=\zeta(s, a)-a^s$. For more details, one can consult \cite{Apostol, Hardy-Wright, Hurwitz}. 

A classical problem which goes back to Goldbach and Bernoulli asserts that

$$\displaystyle\sum_{\omega\in S}(\omega-1)^{-1}=1,$$
\vspace{0.3cm}
where $\displaystyle S=\{n^k: n, k\in\mathbb{Z}_{\geq 0}-\{1\}\}$. In terms of the Riemann zeta function $\zeta(s)$, the above problem reads as,

$$\displaystyle\sum_{n=2}^{\infty}(\zeta(n)-1)=1.$$

Also, there are other representations for $\log 2$ and $\gamma$ (Euler-Mascheroni constant) such as,

$$\displaystyle \sum_{n=1}^{\infty}\frac{\zeta(2n)-1}{n}=\log 2,$$
and

$$\displaystyle \sum_{n=2}^{\infty}\frac{\zeta(n)-1}{n}=1-\gamma.$$
\vspace{0.3cm}

For instance, one way to generate rational zeta series involving $\zeta(2n)$ is by looking at the cotangent power series formula in the form:

$$\displaystyle\sum_{n=1}^{\infty}\zeta(2n)x^{2n}=\frac{1}{2}(1-\pi x\cot(\pi x)), |x|<1.$$

\vspace{0.3cm}

Dividing by $x$ and integrating once, we have

$$\displaystyle \sum_{n=1}^{\infty}\frac{\zeta(2n)}{n}x^{2n}=\log\left(\frac{\pi x}{\sin(\pi x)}\right), |x|<1.$$

\vspace{0.3cm}

For $x=\frac{1}{2}$ and $x=\frac{1}{4}$ in the above formulas, we obtain the following representations:

\begin{equation}
\sum_{n=1}^{\infty}\frac{\zeta(2n)}{2^{2n}}=\frac{1}{2}    
\end{equation}

\begin{equation}
\sum_{n=1}^{\infty}\frac{\zeta(2n)}{2^{4n}}=\frac{4-\pi}{8}    
\end{equation}

\begin{equation}
\sum_{n=1}^{\infty}\frac{\zeta(2n)}{n2^{2n}}=\log\pi-\log 2    
\end{equation}

\begin{equation}
\sum_{n=1}^{\infty}\frac{\zeta(2n)}{n2^{4n}}=\log\pi-\frac{3}{2}\log 2    
\end{equation}

Moreover, integrating from $0$ to $\frac{1}{2}$ the last power series equality, we derive

\begin{equation}
\sum_{n=1}^{\infty}\frac{\zeta(2n)}{n(2n+1)2^{2n}}=\log\pi-1    
\end{equation}
which can be found in \cite{Tyler-Chernoff}.

This type of rational zeta series and many others are treated in \cite{Borwein-Bradley-Crandall}. In \cite{Lupu-Orr} there are given exact formulas for the following rational zeta series

\begin{equation}
\displaystyle\sum_{n=1}^{\infty}\frac{\zeta(2n)}{n4^n}\binom{2n}{m}
\end{equation}

and

\begin{equation}
\displaystyle\sum_{n=1}^{\infty}\frac{\zeta(2n)}{n16^n}\binom{2n}{m}
\end{equation}

\vspace{0.3cm}

in terms of zeta values. In this note, we give an exact formula for a more general rational zeta series which encompasses the two series above. 

The main result of this note is the following

\bt Let $\zeta(s)$ and $\zeta(s; r)$ be the Riemann and Hurwitz zeta functions. For $a\in (0, 1)$ we have

$$\displaystyle \sum_{n=1}^{\infty}\cfrac{a^{2n}\zeta(2n)}{n}\dbinom{2n}{m}=\cfrac{a^m}{m}\Bigg((-1)^m\zeta(m; a)+\zeta(m; 1-a)\Bigg)+\cfrac{(-1)^{m-1}}{m}.$$
\et
\vspace{0.2cm}

The main idea of the proof is a combination of expressing the rational zeta series from the left-hand side as the $n$th derivative of the cotangent function (a polynomial $P_{n}(\cot \pi x)$) and a surprising result of Hoffman \cite{Hoffman} which relates this polynomial $P_{n}(\cot x)$ in terms of Hurwitz zeta values. 

\medskip

\textbf{Acknowledgement.} The second
author was supported by the project “Group schemes, root systems, and related representations”
founded by the European Union - NextGenerationEU through Romania’s National Recovery and
Resilience Plan (PNRR) call no. PNRR-III-C9-2023- I8, Project CF159/31.07.2023, and coordinated by the Ministry of Research, Innovation and Digitalization (MCID) of Romania.

\section{The proof of Theorem 1.1}

Before we dive into the proof of the main result, let us recall a result of Hoffman \cite{Hoffman} which will be an essential ingredient for our purpose.

\bl[M.E. Hoffman, 1995] For real $0<a<1$ and integer $n\geq 0$,

(a) $$\displaystyle\sum_{k=0}^{\infty}\left[\frac{1}{(k+a)^{n+1}}+\frac{(-1)^{n+1}}{(k+1-a)^{n+1}}\right]=\frac{\pi^{n+1}}{n!}P_{n}(\cot a\pi).$$
and

(b) $$\displaystyle\sum_{k=0}^{\infty}\frac{(-1)^k}{(k+a)^{n+1}}+(-1)^n\sum_{k=0}^{\infty}\frac{(-1)^k}{(k+1-a)^{n+1}}=\frac{\pi^{n+1}}{n!}\csc a\pi P_{n}(\cot a\pi)$$

\el

For a function $\alpha$ we will denote by $\alpha^{(k)}$ its $k$th derivative. Our purpose is to explore sums of the type 
$$g(a)=\sum_{n=1}^{\infty}\cfrac{a^{2n}\zeta(2n)}{n}\dbinom{2n}{m}$$
where $a\in (0, 1)$ is a general real parameter.

We start with the following well known cotangent power series expansion written in the form

$$\pi\cot(\pi x)=\cfrac{1}{x}-2\sum_{k=1}^{\infty} \zeta(2k) x^{2k-1}.$$
To connect the above with this expansion, note  that 
$\cfrac{1}{n}\dbinom{2n}{m}=\cfrac{2\cdot(2n-1)\ldots (2n-m+1)}{m!}$.

Thus $g(a)$ can be computed as $$\cfrac{a^{m}}{m!} h^{(m-1)}(a)$$ where $h(x)=-\pi\cot(\pi x)+\cfrac{1}{x}$.

In \cite{Hoffman}, Hoffman computes the $n$th derivative of $-\cot x $ in terms of a polynomial $P_n$. More precisely
$$\displaystyle\cfrac{d^n}{dx^n} \cot x= (-1)^n P_n(\cot x).$$
For our problem we have $\cfrac{d^n}{dx^n} \cot(\pi x)= (-1)^n \pi^n P_n(\cot(\pi x))$. Moreover for $0<a<1$  by Lemma 2.1 (part (a)), we have

$$\sum_{k=0}^{\infty}\left[\cfrac{1}{(k+a)^{n+1}}+\cfrac{(-1)^{n+1}}{(k+1-a)^{n+1}}\right]=
\cfrac{\pi^{n+1}}{(n+1)!}P_n(\cot(\pi a)).$$
Putting everything together we get the following expression for $g(a)$ , which is valid for $0<a<1$
$$g(a)=\cfrac{a^m}{m}\sum_{k=0}^{\infty}\left[\cfrac{(-1)^{m}}{(k+a)^{m}}+\cfrac{1}{(k+1-a)^{m}}\right]+\cfrac{(-1)^{m-1}}{m}.$$
Note that we can compactly write the above in terms of Hurwitz zeta function, 
$$g(a)=\cfrac{a^m}{m}\Bigg((-1)^m\zeta(m,a)+\zeta(m,1-a)\Bigg)+\cfrac{(-1)^{m-1}}{m}.$$$\qed$

\textbf{Remark.} Using Wikipedia for the Hurwitz zeta functions \cite{Hurwitz} one can express for $0<p<q$ and $\gcd(p,q)=1$
$$\zeta\left(s,\frac{p}{q}\right)=\cfrac{q^s}{\varphi(q)}\sum_{\chi}\overline{\chi}(p)L(s,\chi)$$
where the sum runs over all Dirichlet characters mod $q$.

Thus one can express concretely $g(\frac{p}{q})$ in terms of Dirichlet's $L$-functions. Dirichlet $L$-functions \cite{Apostol, Hardy-Wright} are defined as follows. First, consider $\chi$ to be a homomorphism from the units of $\mathbb{Z}/k\mathbb{Z}$ to $\mathbb{C}^{*}$. Now, we can extend $\chi$ to a function on $\mathbb{Z}$ called Dirichlet character modulo $q$ as follows

$$\displaystyle \chi(n)=  \left\{
\begin{array}{ll}
      \displaystyle\chi(q\mathbb{Z}+n)& \operatorname{gcd}(n, k)=1, \\ \\
      \displaystyle 0  &\operatorname{otherwise}.\\
\end{array} 
\right.$$

Then for $s\in\mathbb{C}$, the Dirichlet $L$-series corresponding to the character $\chi$ is given by
$$\displaystyle L(s, \chi)=\sum_{n=1}^{\infty}\frac{\chi(n)}{n^s}=\prod_{p}\left(1-\frac{\chi(p)}{p^s}\right)^{-1}.$$

\vspace{0.2cm}

As corollaries of the main theorem (Theorem 1.1), we have the following representations from \cite{Lupu-Orr}.

\begin{co}

\begin{equation}
\displaystyle\sum_{n=1}^{\infty}\frac{\zeta(2n)}{n4^n}\binom{2n}{m}=  \left\{
\begin{array}{ll}
      \displaystyle\frac{1}{m} & m \operatorname{odd}, \\ \\
      \displaystyle\frac{1}{m}\left(2\zeta(m)\left(1-\frac{1}{2^m}\right)-1\right)  & m \operatorname{even}.\\
\end{array} 
\right.
\end{equation}

\end{co}

and

\begin{co}
We have the following series representation

\begin{equation}
\sum_{n=1}^{\infty}\frac{\zeta(2n)}{n16^n}\binom{2n}{m}=\left\{
\begin{array}{ll}
      \displaystyle \frac{1}{m}\left(1-\beta(m)\right) & m \operatorname{odd}, \\ \\
      \displaystyle \frac{1}{m}\left(\zeta(m)\left(1-\frac{1}{2^m}\right)-1\right)  & m \operatorname{even},\\
\end{array} 
\right.
\end{equation}
where $\displaystyle\beta(s)=\sum_{n=0}^{\infty}\frac{(-1)^n}{(2n+1)^s}$ is the Dirichlet's beta function. 
\bigskip
\end{co}

 As it has been showed in \cite{Lupu-Orr}, many well-known rational zeta series can be obtained from the last two corollaries which are similar with $(1), (2), (3), (4), (5)$. We refer the readers to \cite{Lupu-Orr} for more details.

\vspace*{3mm}
\begin{flushright}
\begin{minipage}{148mm}\sc\footnotesize

Beijing Institute of Mathematical Sciences and Applications (BIMSA),
Yau Mathematical Sciences Center (YMSC), Tsinghua University, Beijing, People's Republic of China\\
{\it E--mail address}: {\tt lupucezar@gmail.com, lupucezar@bimsa.cn} \vspace*{3mm}
\end{minipage}

\begin{minipage}{148mm}\sc\footnotesize

Simion Stoilow Institute of Mathematics of the Romanian Academy, Bucharest, Romania\\
{\it E--mail address}: {\tt vmatei@imar.ro} \vspace*{3mm}
\end{minipage}

\end{flushright}


\begin{thebibliography}{199}
\markboth{Bibliography}{Bibliography}


\bibitem{Apostol}
T. Apostol, {\em Introduction to Analytic Number Theory}, 3rd ed., Springer-Verlag, New York, 1986.

\bibitem{Borwein-Bradley-Crandall}
Jonathan M. Borwein, David M. Bradley, Richard E. Crandall, Computational Strategies for the Riemann Zeta Function, {\em J. Comput. Appl. Math.} \textbf{121} (2000), 247–-296.

\bibitem{Hardy-Wright}
G.H. Hardy, E.M. Wright, {\em An Introduction to the Theory of Numbers}, 4th edn. Oxford University Press, London (1960).

\bibitem{Hoffman}
M. E. Hoffman, Derivative polynomials for tangent and secant, {\em Amer. Math. Monthly} \textbf{102} (1995), 23--30.

\bibitem{Lupu-Orr}
C. Lupu, D. Orr, Series representations for the Apery constant $\zeta(3)$ involving the values $\zeta (2n)$, {\em Ramanujan J.} \textbf{45} (2019), 477--494. 

\bibitem{Tyler-Chernoff}
D. Tyler, P. R. Chernoff, An old sum reappears-Elementary problem 3103, {\em Amer. Math.
Monthly} \textbf{92} (1985), 507.

\bibitem{Hurwitz}
Wikipedia-Hurwitz zeta function, $\mathtt{https://en.wikipedia.org/wiki/Hurwitz\_zeta\_function}$


\end{thebibliography}
\end{document}